\documentclass[conference]{IEEEtran}

\usepackage{algorithm}
\usepackage{algorithmic}
\usepackage{amsmath,amssymb,amsfonts}
\usepackage{mathrsfs}
\usepackage{graphicx}
\usepackage{textcomp}
\usepackage{xcolor}
\usepackage{multirow}
\usepackage{makecell}

\usepackage{nomencl,etoolbox,ragged2e,siunitx,mathtools}

\begin{document}

\title{Quantum Hamiltonian Descent based Augmented Lagrangian Method for Constrained Nonconvex Nonlinear Optimization}

\author{\IEEEauthorblockN{Mingze Li$^{1,2}$,~\textit{Student Member,~IEEE,} Lei Fan$^{1,2}$,~\textit{Senior Member,~IEEE,} and Zhu Han$^{1}$~\textit{Fellow,~IEEE}}
\IEEEauthorblockA{
\textit{Dept. of Electrical and Computer Engineering$^{1}$ and Dept. of Engineering Technology$^{2}$} \\
\textit{University of Houston, TX, USA}\\
\textit{Emails:mli44@cougarnet.uh.edu;
lfan8@central.uh.edu; zhan2@central.uh.edu}
}
}

\maketitle

\begin{abstract}
Nonlinear programming (NLP) plays a critical role in domains such as power energy systems, chemical engineering, communication networks, and financial engineering. However, solving large-scale, nonconvex NLP problems remains a significant challenge due to the complexity of the solution landscape and the presence of nonlinear nonconvex constraints. In this paper, we develop a Quantum Hamiltonian Descent based Augmented Lagrange Method (QHD-ALM) framework to address large-scale, constrained nonconvex NLP problems. The augmented Lagrange method (ALM) can convert a constrained NLP to an unconstrained NLP, which can be solved by using Quantum Hamiltonian Descent (QHD). To run the QHD on a classical machine, we propose to use the Simulated Bifurcation algorithm as the engine to simulate the dynamic process.  We apply our algorithm to a Power-to-Hydrogen System, and the simulation results verify the effectiveness of our algorithm.
\end{abstract}


\section{Introduction}
Given its generality and modeling flexibility, nonlinear programming (NLP) has become a core tool across a wide spectrum of application domains. In power and energy systems, it plays a central role in optimal power flow, unit commitment, and the integration of renewable resources \cite{10906982}. A variety of solvers like IPOPT, Gurobi, and CPLEX have been developed to tackle specific instances of NLP problems.  However, these classical algorithms can only ensure local convergence for nonconvex problems. Many challenges remain in ensuring scalability, robustness, and the ability to escape local minima in large-scale, nonconvex settings.

To address the limitations of classical solvers in navigating nonconvex and high-dimensional optimization landscapes, recent research has turned to quantum computing approaches. One promising direction is quantum annealing, which leverages quantum tunneling to escape local minima and more effectively explore the global structure of the objective landscape \cite{10841900,9771632}. Quantum annealers have been explored for solving various optimization problems, including discrete and combinatorial optimization problems. Extending this concept into continuous domains, the Quantum Hamiltonian Descent (QHD) algorithm has been proposed to solve box-constrained continuous nonlinear optimization problems \cite{leng2023quantumhamiltoniandescent,leng2025quantumhamiltoniandescentnonsmooth}. By exploiting quantum mechanisms such as wavefunction evolution and tunneling-like dynamics, QHD can, in some cases, outperform classical solvers, especially for nonconvex problems where traditional gradient-based methods tend to get stuck in poor local minima. In addition, recent advancements in Hamiltonian embedding techniques have significantly enhanced the feasibility of implementing quantum simulations on near-term quantum devices \cite{leng2024expandinghardwareefficientlymanipulablehilbert}. QHD has been implemented in a software package called the Quantum Hamiltonian Descent Optimizer (QHDOPT) \cite{kushnir2024qhdopt}. This package can interface with various quantum hardware backends to solve box-constrained nonlinear nonconvex optimization problems. QHDOPT can encode the original problem into an Ising model, which can then be solved by quantum annealing machines. 

Despite the theoretical advantages of the quantum adiabatic process, current quantum hardware suffers from significant practical limitations. The number of available qubits remains limited, and noise, decoherence, and the lack of full programmability constrain the size and complexity of problems that can be realistically solved. To overcome these hardware barriers while retaining the benefits of quantum dynamics, researchers have proposed classical algorithms that mimic quantum behaviors. One approach is the simulated bifurcation (SB) algorithm, which emulates the adiabatic evolution of a Hamiltonian system using a network of classical nonlinear oscillators \cite{doi:10.1126/sciadv.aav2372}. SB can be interpreted as a form of Hamiltonian descent, capturing essential features of quantum annealing—such as bifurcation-driven transitions to low-energy states—without requiring quantum hardware. Compared with classical simulated annealing, SB exhibits improved performance on large-scale combinatorial problems by efficiently exploring complex energy landscapes through collective dynamics. Its inherently parallel structure and ability to exploit physical analogies make SB particularly well-suited for integration into quantum optimization solvers like QHDOPT. 

While QHDOPT is well-suited for solving unconstrained or box-constrained nonlinear problems, many real-world applications involve more general equality or inequality constraints that fall outside this scope. In this paper, we integrate QHDOPT into the Augmented Lagrangian Method (ALM) framework. In the ALM framework, a constrained optimization problem is reformulated as a box-constrained problem by introducing Lagrange multipliers and penalty terms into the objective function. This allows the original constraints to be enforced indirectly, as violations are penalized during optimization. By progressively updating the Lagrange multipliers and increasing the penalty parameters, the algorithm converges toward a solution that satisfies both optimality and feasibility. The newly proposed QHD based ALM framework can tackle a broader class of constrained nonlinear programming problems while maintaining its ability to leverage quantum-inspired dynamics for global optimization.

In Section~\ref{sec:2}, we present the general formulation of NLP problems and discuss the challenges associated with solving nonconvex, large-scale instances. Section~\ref{sec:3} introduces the quantum optimization framework QHDOPT, which models continuous box-constrained NLP problems as quantum dynamical systems. Section~\ref{sec:4} incorporates the SB algorithm, a classical heuristic that enables QHDOPT to run on classical machines. In Section~\ref{sec:5}, we introduce the QHD-ALM framework. In Section~\ref{sec:6}, we apply the proposed QHD-ALM to solve a hydrogen energy production management problem.

\section{Introduction to Nonlinear Programming (NLP)}\label{sec:2}

Nonlinear programming is a mathematical optimization framework used to solve problems where the objective function and constraints involve nonlinear relationships. Unlike linear programming, which assumes a linear structure, NLP allows for a more flexible representation of practical industrial systems. However, the inherent complexity of nonlinear models often leads to computational challenges, especially when dealing with nonconvexity, where multiple local optima exist, making it difficult to find the global solution.

A general NLP problem is formulated as follows:
\begin{align}
    \min_{\textbf{x}} &\quad f(\textbf{x}) \label{eq:objective} \\
    \text{s.t.} &\quad g_i(\textbf{x}) = 0, \quad i \in \mathcal{I}, \label{eq:equality} \\
    &\quad h_j(\textbf{x}) \leq 0, \quad j \in \mathcal{J}, \label{eq:inequality}
\end{align}
where \( \textbf{x} \) is a vector of continuous decision variables, \( f(\textbf{x}) \) is a nonlinear objective function, and \( g_i(\textbf{x}) \) and \( h_j(\textbf{x}) \) represent nonlinear equality and inequality constraints, respectively. The nonlinearity in these functions may arise from physical relationships, efficiency curves, or interactions between different components in the system. If the objective function or constraints are nonconvex, the problem becomes significantly harder to solve since classical optimization techniques are no longer guaranteed to find the global optimum. To overcome these challenges, advanced computational methods are needed. Quantum algorithms, such as QHD, exploit quantum mechanical properties like tunneling to escape local minima and explore a broader solution space. 

Given the limitations of traditional solvers in handling large-scale non-convex NLP problems, this paper introduces an advanced algorithmic framework designed to enhance the solution process. By integrating QHD and ALM, we aim to improve the efficiency and reliability in solving NLP problems with nonlinear and non-convex constraints. 

\section{Introduction to QHDOPT} \label{sec:3}
QHDOPT is an open-source optimization solver that implements QHD to solve continuous NLP problems with box constraints. The key idea behind QHD is to model the optimization process as the evolution of a quantum wavefunction governed by a time-dependent Hamiltonian. Unlike traditional optimization methods that rely on gradient-based search, QHD introduces a dynamical system where both kinetic and potential energy components play a role in guiding the optimization trajectory. This quantum formulation enables QHD to leverage quantum tunneling, providing an efficient means of escaping local minima in highly nonconvex problems \cite{leng2023quantumhamiltoniandescent,leng2025quantumhamiltoniandescentnonsmooth,kushnir2024qhdopt}.

QHDOPT is designed to handle optimization problems of the form in \cite{kushnir2024qhdopt}:
\begin{equation}
    \min_{x \in \mathbb{R}^n} f(x) = \sum_{i=1}^{n} g_i(x_i) + \sum_{j=1}^{m} p_j(x_{k_j}) q_j(x_{\ell_j}),
\end{equation}
subject to box constraints:
\begin{equation}
    L_i \leq x_i \leq U_i, \quad \forall i=1,\dots,n.
\end{equation}
Here, the objective function consists of both univariate terms \( g_i(x_i) \) and bivariate interaction terms \( p_j(x_{k_j}) q_j(x_{\ell_j}) \). To encode the objective function into the Hamiltonian, QHD constructs the Hamiltonian as follows:
\begin{equation}
H(t)=e^{\varphi t}\left(-\frac{1}{2} \Delta\right)+e^{\chi t}f(x),
\end{equation}
where \( x \) represents the optimization variables, \( \Delta \) denotes the Laplacian operator. The time-dependent scaling factors \( e^{\varphi t} \) and \( e^{\lambda t} \) govern the total energy distribution of the quantum system. This mechanism allows QHD to function as a quantum-enhanced optimization process, dynamically adjusting its exploration-exploitation balance while avoiding the pitfalls of local minima.

To implement QHD in the discrete time domain, QHDOPT applies spatial discretization to represent the continuous wavefunction over a finite grid in \( \mathbb{R}^n \). Each variable domain is discretized into \( N \) grid points, and differential operators such as the Laplacian are approximated using finite-difference schemes. The objective function is encoded as a potential energy operator acting on the discretized space. The resulting QHD Hamiltonian is a large Hermitian matrix that governs quantum evolution over a finite-dimensional state space. 
\begin{equation}
\hat{H}(t)=e^{\varphi_t}\left(-\frac{1}{2} L_d\right)+e^{\chi_t} F_d,
\end{equation}
where
\begin{equation}
\begin{gathered}
L_d=\sum_{i=1}^n I \otimes \cdots \otimes L \otimes D\left(g_i\right)\otimes \ldots I, \\
F_d=\sum_{i=1}^n I \otimes \cdots \otimes D\left(p_j\right) \otimes \sum_{j=1}^m I \otimes \cdots \otimes D\left(q_j\right)\otimes \ldots I .
\end{gathered}
\end{equation}
Here, $I$ is the $N$-dimensional identity operator, $L$ and $D(*)$ denote $N$-dimensional matrices.

To further accelerate the simulation of this Hamiltonian on quantum hardware or classical emulators, QHDOPT adopts a Hamiltonian embedding strategy that converts the high-dimensional QHD Hamiltonian into an Ising-type representation suitable for quantum machines \cite{leng2024expandinghardwareefficientlymanipulablehilbert}.

The discretized QHD Hamiltonian is a Hermitian matrix with an explicit tensor product decomposition structure. This structure enables the application of Hamiltonian embedding techniques to construct a surrogate Hamiltonian. By mapping each discretized variable to a binary encoding and converting binary variables to spin variables \( s_i \in \{-1, 1\} \), the surrogate Hamiltonian can be embedded into an Ising model which can be solved by quantum solvers.

\section{Simulated Bifurcation Algorithm for Hamiltonian Simulation in QHDOPT} \label{sec:4}

In earlier implementations of QHDOPT, quantum annealers were used to solve Ising models that arise within the QHD framework. However, current quantum hardware remains constrained by noise, limited qubit counts, and a lack of fault tolerance, posing practical challenges for scalable deployment. To overcome these limitations while preserving the quantum-inspired philosophy of Hamiltonian descent, we replace the quantum annealing step with the Simulated Bifurcation (SB) algorithm—a purely classical alternative.

The SB algorithm is a heuristic method inspired by quantum adiabatic evolution but implemented entirely in a classical framework. It efficiently solves the Ising model by leveraging a network of coupled Duffing oscillators and their bifurcation dynamics. While quantum annealing relies on tunneling to overcome energy barriers, SB exploits the classical bifurcation mechanism of nonlinear oscillators to traverse complex energy landscapes. This allows highly parallel and scalable execution on modern classical hardware such as GPUs. In QHDOPT,  we use SB to replace quantum annealing to perform Hamiltonian descent efficiently, retaining the benefits of quantum-inspired optimization with enhanced computational accessibility.

SB encodes the objective into a Hamiltonian system of nonlinear oscillators, evolving the system from a symmetric initial state to a bifurcated state representing an optimal solution. This approach is rooted in quantum adiabatic optimization using Kerr-nonlinear parametric oscillators (KPOs), where each binary variable in the Ising model corresponds to a KPO. The SB algorithm is governed by the simplified Hamiltonian:
\begin{align}
H_q(t) = &\hbar \sum_{i=1}^{N} \left[ \frac{K}{2} a_i^{\dagger 2} a_i^2 - \frac{p(t)}{2}(a_i^{\dagger 2} + a_i^2) + \Delta_i a_i^{\dagger} a_i \right] \notag\\
&-\hbar \xi_0 \sum_{i=1}^{N} \sum_{j=1}^{N} I_{ij} a_i^{\dagger} a_j,
\end{align}
where $\hbar$ is the reduced Planck constant, \( a_i^\dagger \) and \( a_i \) are the creation and annihilation operators of the \( i \)th oscillator, \( K \) is the Kerr coefficient, \( \Delta_i \) is the detuning, \( p(t) \) is the time-dependent pumping amplitude, and \( \xi_0 \) controls the coupling strength between oscillators through the Ising interaction matrix \( I_{ij} \).

As the system evolves, the amplitude \( p(t) \) is increased slowly from zero, and each KPO transitions from the vacuum state to a coherent state with either positive or negative amplitude. According to the quantum adiabatic theorem, the system remains in its instantaneous ground state, and the signs of the final coherent amplitudes encode the Ising spins. This mechanism provides a quantum dynamical route to solving Ising problems.

To implement this on classical hardware, a classical approximation is taken by replacing quantum operators with continuous-valued expectation variables. Specifically, each operator \( a_i \) is approximated by a complex amplitude \( x_i + i y_i \), where \( x_i \) and \( y_i \) are interpreted as canonical position and momentum variables. This leads to the classical mechanical Hamiltonian used in the SB algorithm. Through this approximation, the SB method retains the quantum-inspired bifurcation dynamics while enabling efficient simulation on classical architectures.

The SB algorithm evolves a system of coupled oscillators governed by the following simplified Hamiltonian:
\begin{align}
H_{\text{SB}}(\mathbf{x}, \mathbf{y}, t) = &\sum_{i=1}^{N} \frac{\Delta}{2} y_i^2 + \sum_{i=1}^{N} \left[ \frac{K}{4} x_i^4 + \frac{\Delta - p(t)}{2} x_i^2 \right] \notag\\
&- \frac{\xi_0}{2} \sum_{i=1}^{N} \sum_{j=1}^{N} I_{ij} x_i x_j,
\end{align}
where \( \mathbf{x} \) and \( \mathbf{y} \) denote position and momentum variables, respectively. The matrix \( I \) encodes problem-specific interactions, and \( p(t) \) is a bifurcation-driving parameter that increases during the simulation. All the detunings have been assumed to be the same value $\Delta$. The system evolves according to Hamilton’s equations \cite{doi:10.1126/sciadv.aav2372}:
\begin{align}
\dot{x}_i &= \Delta y_i, \\
\dot{y}_i &= -\left[ K x_i^3 - (p(t) - \Delta) x_i + \xi_0 \sum_{j=1}^{N} I_{ij} x_j \right].
\end{align}

These equations are numerically integrated using the explicit symplectic Euler method, which ensures numerical stability and allows for large time steps. All oscillator states are initialized near zero, and as \( p(t) \) increases, the system undergoes bifurcation. At the final time step, the spin configuration is determined by the sign of each oscillator position.

The SB algorithm efficiently explores multiple minima
and escapes poor local optima. Its simplicity allows for fast
computation, with parallelization enabled by its structure.
By replacing quantum annealing with SB, QHDOPT gains
scalability, stability, and flexibility, enabling the solution of large-scale nonlinear optimization problems that would be
difficult for current quantum devices.

\section{Integration of QHD with Augmented Lagrangian Method} \label{sec:5}
In this section, we will introduce how to integrate the QHDOPT into the ALM framework.

\subsection{General Formulation of ALM}
We consider a general NLP problem with both equality and inequality constraints, written in the standard form as shown in \eqref{eq:objective}, \eqref{eq:equality} and \eqref{eq:inequality}. To handle the inequality constraints, we introduce nonnegative slack variables \( s_j \geq 0 \). This reformulates the original problem into an equality-constrained problem, where all inequality constraints are rewritten as
\begin{align}
    g_i(x) &= 0, \quad i \in \mathcal{I}, \\
    h_j(x) + s_j &= 0, \quad s_j \geq 0, \quad j \in \mathcal{J}.
\end{align}

The augmented Lagrangian function is then defined as:
\begin{equation}
\begin{aligned}
    \mathcal{L}_A(x, s, \lambda, \mu, \rho) &= f(x) 
    + \sum_{i \in \mathcal{I}} \lambda_i g_i(x) 
    + \sum_{j \in \mathcal{J}} \mu_j (h_j(x) + s_j) \\
    &\quad +  \sum_{i \in \mathcal{I}} \frac{\rho_i}{2}g_i(x)^2 
    +  \sum_{j \in \mathcal{J}} \frac{\rho_j}{2}(h_j(x) + s_j)^2,
\end{aligned}
\label{eq:alm-formulation}
\end{equation}
where \( \lambda_i \) and \( \mu_j \) are the Lagrange multipliers for the equality and transformed inequality constraints, and \( \rho_{i}, \rho_{j} > 0 \) are penalty parameters.

The ALM algorithm proceeds by solving a sequence of unconstrained problems:
\[
(x^{(k)}, s^{(k)}) = \arg\min_{x, s} \mathcal{L}_A(x, s, \lambda^{(k)}, \mu^{(k)}, \rho^{(k)}),
\]
and then updating the multipliers after each iteration using \cite{5570998}:
\begin{align}
    \lambda_i^{(k+1)} &= \lambda_i^{(k)} + \rho_i^{(k)} g_i(x^{(k)}), \quad i \in \mathcal{I}, \\
    \mu_j^{(k+1)} &= \mu_j^{(k)} + \rho_j^{(k)} (h_j(x^{(k)}) + s_j^{(k)}), \quad j \in \mathcal{J}.
\end{align}

The penalty parameter \( \rho^{(k)} \) is optionally increased over iterations to improve convergence.

This framework enables the use of box-constrained or unconstrained solvers, such as QHDOPT or IPOPT, to efficiently solve the reformulated problem. In particular, for highly nonlinear or nonconvex problems, the use of QHD in the ALM framework can improve convergence and solution quality while maintaining the feasibility of the original constraints.

\subsection{Advantages of Using QHD-ALM}
\begin{figure}[t]
    \centering
\includegraphics[scale=0.35]{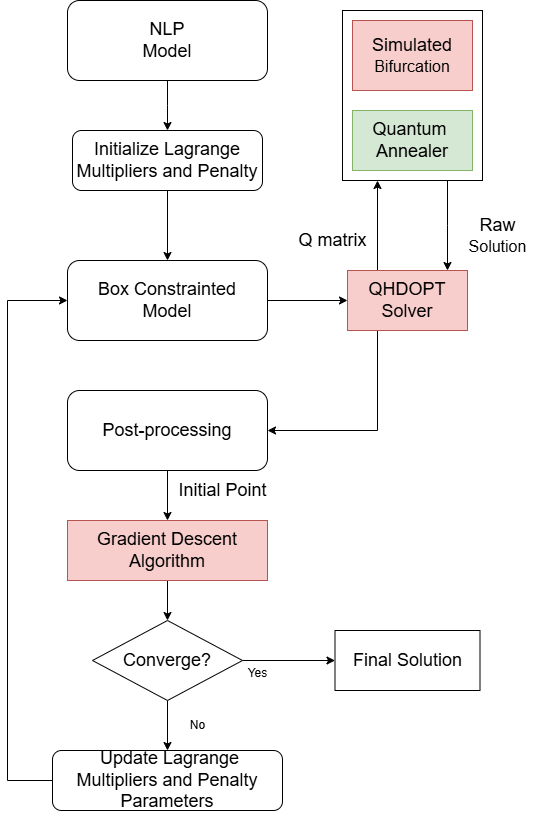}
    \caption{Flow Chart of QHD-ALM Method}
    \label{fig:flow}
\end{figure}
Fig.~\ref{fig:flow} and Algorithm~\ref{alg:qhdopt-alm} illustrate the QHD-ALM framework. To enable compatibility with quantum solvers, the method first initializes Lagrange multipliers and penalty parameters, preparing the problem for transformation into a box-constrained form via ALM. This box-constrained model is then passed to the QHDOPT solver, which generates a coefficient matrix representing the discretized optimization landscape. To solve this binary-encoded representation, QHDOPT can interface with a quantum annealer or a classical algorithm such as the simulated bifurcation method. These solvers return raw candidate solutions, which are then refined in a post-processing step to improve feasibility and extract a suitable initial point for further local optimization. This initial point is used to warm-start a classical solver, such as IPOPT, which performs fine-grained numerical optimization on the original constrained problem. If the solution has not yet converged to an acceptable threshold, the framework updates the Lagrange multipliers and penalty parameters, and repeats the process iteratively. This hybrid loop continues until convergence criteria are met, at which point a high-quality final solution is returned.
\begin{algorithm}
    \caption{QHD-ALM Framework }
    \label{alg:qhdopt-alm}
    \begin{algorithmic}[1]
        \STATE \textbf{Input:} Nonlinear programming model with constraints
        
        \STATE \textbf{Step 1: Initialize}
        \STATE Set initial Lagrange multipliers \( \lambda^{(0)} \), \( \mu^{(0)} \), and penalty \( \rho^{(0)} \)
        \STATE Set iteration counter \( k = 0 \)

        \REPEAT
            \STATE \textbf{Step 2: Box-constrained Model Reformulation}
            \STATE Construct the Augmented Lagrangian function \eqref{eq:alm-formulation}
            \STATE \textbf{Step 3: Solve with QHDOPT}
            \STATE Formulate Q matrix from the box-constrained model
            \STATE Use QHDOPT with either Simulated Bifurcation or Quantum Annealer to obtain a raw solution
            
            \STATE \textbf{Step 4: Post-processing}
            \STATE Map the raw solution to original feasible space
            \STATE Use it as an initial point

            \STATE \textbf{Step 5: Refinement with Gradient Descent Algorithm}
            \STATE Utilize gradient descent algorithm to solve the unconstrained NLP from an initial point

            \STATE \textbf{Step 6: Check Convergence}
            \IF{convergence criterion is met}
                \STATE \textbf{Output:} Final solution
                \STATE \textbf{Break}
            \ELSE
                \STATE \textbf{Step 7: Update Parameters}
                \STATE Update \( \lambda^{(k+1)}, \mu^{(k+1)}\)
                \STATE Increase penalty \( \rho^{(k+1)} = \gamma \rho^{(k)} \), where \(\gamma > 1\)
                \STATE \( k \leftarrow k + 1 \)
            \ENDIF
        \UNTIL{convergence is achieved}
    \end{algorithmic}
\end{algorithm}

\section{Case Study}\label{sec:6}
In this section, we apply our proposed QHD-ALM to solve a hydrogen production problem.
\subsection{Mathematical Model}
The objective of the optimization problem is to minimize the operational cost of hydrogen production. In the objective function \eqref{obj}, the total electricity cost over the planning horizon is penalized, while the value of stored hydrogen at the final time step is rewarded. Here, \( C^{\mathrm{hyo}} \) is the market value of hydrogen, and \( C_t^{\mathrm{power}} \) represents the electricity price at time step \( t \).
\begin{subequations}
\begin{align}
\max\ & C^{\mathrm{hyo}} (s_N - s_0) - \sum_{t=0}^{N-1} C_t^{\mathrm{power}} p_t^{\mathrm{buy}}, \label{obj} \\
\text{s.t. } & p_t^{\mathrm{buy}} + p_t^{\mathrm{R}} = m^{\mathrm{AC}} p_t^{{el}} + k^{\mathrm{AC}}, \label{power} \\
& s_{t+1} = s_t + s_t^{{el}} - s_t^{\mathrm{d}}, \label{stor} \\
& S^{\min} \leq s_t \leq S^{\max}, \label{stormax} \\
& 0 \leq p_t^{{el}} \leq P^{\max}, \label{pmax} \\
& p_t^{\mathrm{buy}} \geq 0, \\
& s_t^{{el}} = \Delta t \cdot \frac{p_t^{{el}}}{HHV_{H_2}} \cdot \lambda_{el}, \quad s_t^{{el}} \geq 0. \label{elec}
\end{align}
\end{subequations}

Constraint \eqref{power} captures the power balance for the electrolyzer, where power demand is met by the grid power \( p_t^{\mathrm{buy}} \) and renewable power \( p_t^{\mathrm{R}} \). The coefficients \( m^{\mathrm{AC}} \) and \( k^{\mathrm{AC}} \) are linear conversion parameters reflecting losses or auxiliary loads. The storage dynamics in \eqref{stor} ensure that the hydrogen storage level \( s_t \) evolves based on production \( s_t^{{el}} \) and consumption \( s_t^{\mathrm{d}} \). Constraints \eqref{stormax} and \eqref{pmax} enforce upper and lower bounds on storage and electrolyzer power capacity, respectively. Constraint \eqref{elec} links the produced hydrogen energy \( s_t^{{el}} \) to the electrolyzer’s input power \( p_t^{{el}} \), scaled by the higher heating value of hydrogen \( HHV_{H_2} \), time duration \( \Delta t \), and electrolyzer efficiency \( \lambda_{el} \).

The efficiency \( \lambda_{el} \) is a critical nonlinear factor influenced by the operating point of the electrolyzer, such as input power and temperature. While \eqref{elec} uses a generalized fixed efficiency for simplicity, the subsequent section develops a detailed \textit{dynamic efficiency model} that incorporates nonlinear characteristics. This refined model improves fidelity by capturing variations in \( \lambda_{el} \) as a function of key physical and operational parameters, thereby enhancing the accuracy and realism of the overall optimization framework.

We have investigated the relationship between power and efficiency, illustrating the performance characteristics under different operating conditions. We perform the curve fitting and obtain the fitted voltage and current in \cite{li2025integrated}:
\begin{align}
& \bar{I}_{cell} =  i_1-i_2\cdot \exp(i_3\cdot p_t^{el})+i_4\cdot p_t^{el},\\
& \bar{U}_{cell} = u_1 + u_2 p_t^{el} - u_3 (p_t^{el})^2 + u_4 (p_t^{el})^3 - u_5 (p_t^{el})^4.
\end{align}
Then model the efficiency as:
\begin{align}
 \lambda_{el} \leq &m_1 + m_2 \cdot P^{\max } + \notag \\
&m_3 \cdot \exp\left(m_4 \cdot (100 \cdot p_t^{{el}}/P^{\max })\right), \label{formu:e_1}\\ 
 \lambda_{el} \leq  &\left(n_1 + n_2 \cdot \exp{\left(\frac{n_3 + n_4 T + n_5 T^2}{\bar{I}_{cell}}\right)} \right)/\bar{U}_{cell},\label{formu:e_2}
\end{align}
where $i,u,m,n,T$ are parameters of the dynamic efficiency functions.

To solve the nonlinear constrained hydrogen scheduling problem using box-constrained solvers, we employ the ALM. The ALM incorporates Lagrange multipliers and penalty terms that penalize violations of the constraints, enabling the problem to be treated as an unconstrained optimization task. The key advantage of this approach is that it allows us to apply box-constrained solvers, which are often more efficient and widely available, to complex constrained optimization problems.

In our case, we are dealing with a nonlinear constrained hydrogen scheduling problem, which involves both nonlinearities and multiple constraints. To apply box-constrained solvers, we need to transform these constraints into penalty terms, which is where the QHD algorithm comes into play. The QHD algorithm is specifically designed to solve box-constrained nonconvex problems, making it a suitable choice for our problem. However, to apply QHD effectively, we must first reformulate the problem by transforming the original constraints into penalty terms.

Let \( x = \{p_t^{{el}}, p_t^{\mathrm{buy}}, s_t, \lambda_t\}_{t=0}^N \) denote the decision variables. We define the following constraints in the standard NLP format:
\begin{align}
    g_{1,t}(x) &= s_{t+1} - s_t - \frac{\Delta t}{HHV_{H_2}} p_t^{{el}} \cdot \frac{\lambda_t}{100} + s_t^{\mathrm{d}}. \\
    h_{1,t}(x) &= \lambda_t - \frac{n_1 + n_2 \cdot \exp\left( \frac{n_3 + n_4 T + n_5 T^2}{\bar{I}_{cell}} \right)}{\bar{U}_{cell}},
    \\
    h_{2,t}(x) &= \lambda_t - \left( m_1 + m_2 P^{\max} + m_3 \exp\left( m_4 \cdot \frac{100 p_t^{{el}}}{P^{\max}} \right) \right) .
\end{align}
To convert \( h_i(x) \leq 0 \) into an equality, we introduce slack variables \( sl_{i,t} \geq 0 \). The resulting Augmented Lagrangian function becomes:
\begin{equation}
\begin{aligned}
\mathcal{L}_A =\ 
& - C^{\mathrm{hyo}} (e_N - e_0) + \sum_{t=0}^{N-1} C_t^{\mathrm{power}} p_t^{\mathrm{buy}} \\
&+ \sum_{t=0}^{N-1} \lambda_{1,t} g_{1,t}(x) + \frac{\rho_1}{2} g_{1,t}(x)^2 \\
&+ \sum_{t=0}^{N-1} \mu_t (h_{1,t}(x) + sl_{1,t}) + \frac{\rho_2}{2} (h_{1,t}(x) + sl_{1,t})^2 \\
&+ \sum_{t=0}^{N-1} \mu_t (h_{2,t}(x) + sl_{2,t}) + \frac{\rho_3}{2} (h_{2,t}(x) + sl_{2,t})^2,
\end{aligned}
\end{equation}
where the objective function includes hydrogen production costs, power purchasing costs, and penalty terms for the constraints. The optimization problem is subject to the following bounds:
\begin{align}
    0 &\leq p_t^{{el}} \leq P^{\max}, S^{\min} \leq s_t \leq S^{\max},
    p_t^{\mathrm{buy}} \geq 0, \notag \\
    0 &\leq \lambda_t \leq 100, 
    sl_{i,t} \geq 0.
\end{align}
This formulation can be solved using QHDOPT, which is well-suited for solving nonlinear, non-convex optimization problems, leveraging quantum optimization advantages such as enhanced exploration of complex solution spaces.

\subsection{Simulation Results}
In the case study, we evaluate four optimization methods across problem instances of varying sizes, corresponding to 6, 24, 168, and 720 time slots. These cases reflect increasing problem dimensions and computational complexity. The methods compared are: (1) pure IPOPT, which solves the problem from a single random initial point; (2) IPOPT with multiple initial points, where the best solution is selected from 1,000 random initializations; (3) the classical Augmented Lagrangian Method; and (4) our proposed QHD-ALM algorithm. For each case, we compare the objective values achieved and the computational time required, highlighting both solution quality and efficiency.

\begin{table}[ht]
\centering
\caption{Optimal Objective Value of Different Methods}
\label{tab:2}
\begin{tabular}{c|c|c|c|c}
\hline
Case & Pure-IPOPT (\$)  & \begin{tabular}[c]{@{}c@{}}IPOPT 1k  \\ Samples (\$)\end{tabular} & ALM (\$)& QHD-ALM (\$)  \\ \hline
1    & 6.42            & 892.06                                                & 6.42       & 893.8 \\
2    & 6.33            & 2312.92                                                & 6.33       & 2333.4 \\
3    & -760.23          & 14153.52                                              & 10124.05   & 13877 \\
4    & 3040.31          & 19368.54                                               & 17423.74   & 18840.1        \\ 
\hline
\end{tabular}
\end{table}

\begin{table}[ht]
\centering
\caption{Computation Time of Different Methods}
\label{tab:1}
\begin{tabular}{c|c|c|c|c}
\hline
Case & Pure-IPOPT  & \begin{tabular}[c]{@{}c@{}}IPOPT 1k  \\ Samples \end{tabular} & ALM& QHD-ALM  \\ \hline
1    & 0.089s          & 73s                                                 & 1.31s   & 6.82s \\
2    & 0.289s            & 257s                                                & 3.18s      & 11.28s \\
3    & 1.767s           & 16 min                                                & 61.1s   & 78.2s \\
4    & 3.184s           & 52 min                                                & 351.58s   & 369.38s         \\ 
\hline
\end{tabular}
\end{table}
Table~\ref{tab:2} illustrates the objective values obtained across different methods for various problem sizes. As shown, the IPOPT solver often struggles with nonconvexity, frequently converging to poor local minima when initialized from a single random point. Although running IPOPT 1,000 times with different initializations improves the solution quality, this brute-force approach is computationally intensive and inefficient. 

Table~\ref{tab:1} compares the total runtime of the four optimization methods across increasing problem sizes. As expected, pure IPOPT is the fastest since it performs only a single local search from a random starting point. However, to reliably find high-quality solutions for nonconvex problems, IPOPT must be executed multiple times with different initializations—resulting in significantly higher total computation time. In contrast, QHD-ALM achieves comparable or better objective values with much less computational cost than IPOPT with 1,000 restarts, offering a more time-efficient strategy for global exploration.

\section{Conclusion}
In this paper, we proposed a novel hybrid optimization framework named QHD-ALM, which integrates QHD into ALM. The proposed algorithm improves the convergence of the classical ALM framework. The application of our proposed algorithm to the hydrogen production management verified the efficiency and optimality of our algorithm.

\end{document}